\documentclass[11pt,a4paper]{article}
\usepackage[T2A]{fontenc}
\usepackage[cp1251]{inputenc}
\usepackage[english]{babel}
\usepackage{epsfig,amssymb,epic,eepic,color, graphicx}
\usepackage{fancyhdr}
\usepackage{amsthm}
\usepackage{amsmath}
\usepackage{indentfirst}
\usepackage{float}

\usepackage{pb-diagram}

\usepackage{authblk}

\usepackage{setspace}
\usepackage{float}
\usepackage{fancyhdr}
\usepackage{geometry}
\newtheorem{theorem}{Theorem}
\newtheorem{corollary}{Corollary}

\newtheorem{lemma}{Lemma}[section]

\theoremstyle{definition}

\newtheorem{remark}{Remark}

\newenvironment{demo}{{\bf Proof: }}{\hfill $\square$ \medskip}
\newenvironment{demo1}{{\bf Proof of Theorem 1: }}{\hfill $\square$ \medskip}
\newenvironment{demo2}{{\bf Proof of Theorem 2: }}{\hfill $\square$ \medskip}
\newenvironment{demo3}{{\bf Proof of Corollary 1: }}{\hfill $\square$ \medskip}
\renewenvironment{abstract}{\textbf{Abstract.}}{\medskip}
\newtheorem*{ac}{\textbf{Acknowledgments}}{\medskip}
\date{}

\begin{document}

\title{On decomposition of ambient surfaces admitting $A$-diffeomorphisms with nontrivial attractors and repellers}

\author{V. Grines, D. Mints \\ HSE University}

\maketitle

\begin{abstract}
It is well-known that there is a close relationship between the dynamics of diffeomorphisms satisfying the axiom $A$ and the topology of the ambient manifold. In the given article, this statement is considered for the class $\mathbb G(M^2)$ of $A$-diffeomorphisms of closed orientable surfaces such that their non-wandering set consists of $k_f\geq 2$ connected components of one-dimensional basic sets (attractors and repellers). We prove that the ambient surface of every diffeomorphism $f\in \mathbb G(M^2)$ is homeomorphic to the connected sum of $k_f$ closed orientable surfaces and $l_f$ two-dimensional tori such that the genus of each surface is determined by the dynamical properties of appropriating connected component of a basic set and $l_f$ is determined by the number and position of bunches, belonging to  all connected components of basic sets. We also prove that every diffeomorphism from the class $\mathbb G(M^2)$ is $\Omega$-stable but is not structurally stable.
\end{abstract}

\section{Introduction and statement of results}\label{s1}

Let $M^n$ be a closed smooth manifold of dimension $n\ge 1$, $f: M^n\rightarrow M^n$ be a diffeomorphism, and $NW (f)$ be a non-wandering set of $f$.

The diffeomorphism $f$ is said to be $A$-diffeomorphism (diffeomorphism satisfying the axiom $A$) if its non-wandering set $NW(f)$ is hyperbolic and the periodic points are everywhere dense in $NW (f)$. According to \cite{1} (Theorem 6.2), the set $NW (f)$  can be uniquely expressed as the finite union of mutually disjoint subsets $\Omega_1,...,\Omega_k$, called basic sets, each of which is compact, invariant and topologically transitive. A basic set which is a periodic trajectory is called trivial. Otherwise, a basic set is called nontrivial.

According to \cite{14}, \cite{13}, any basic set $\Omega_i$ $(i\in\{1,..., k\})$ of the diffeomorphism $f: M^n\rightarrow M^n$ is uniquely expressed as the finite union of the compact subsets

$$\Omega_i=\Omega_{i_1}\cup...\cup\Omega_{i_q}, \; q\ge 1,$$

called periodic components of the set $\Omega_i$\footnote{R. Bowen \cite{13} called these components $C$-dense. In this paper, following \cite{5}, we call them periodic (by analogy with periodic points of a periodic orbit).}, such that $f^q (\Omega_{i_j})=\Omega_{i_j}, f(\Omega_{i_j})=\Omega_{i_{j+1}}, j\in\{1,..., q\}$ $(\Omega_{i_{q+1}}=\Omega_{i_{1}})$. For every point $x$ belonging to the periodic component $\Omega_{i_j}$, the set $W^s_x\cap\Omega_{i_j}$ $(W^u_x\cap\Omega_{i_j})$ is dense in $\Omega_{i_j}$, where $W^s_x$ $(W^u_x)$ denotes a stable (unstable) manifold of the point $x$.

A basic set $\Omega_i$ is called an attractor (repeller) if it has a closed trapping neighborhood $U_{\Omega_i}\subset M^n$ such that $f(U_{\Omega_i})\subset int$ $U_{\Omega_i}, \bigcap\limits_{k\in\mathbb N} f^k(U_{\Omega_i})=\Omega_i$ $(f^{-1}(U_{\Omega_i})\subset int$ $U_{\Omega_i}, \bigcap\limits_{k\in\mathbb N} f^{-k}(U_{\Omega_i})=\Omega_i)$. In this case, $\Omega_i = \bigcup\limits_{x\in\Omega_i}W^u (x)$ $(\Omega_i = \bigcup\limits_{x\in\Omega_i}W^s (x))$. If $ \dim\Omega_i=\dim W^u (x)$ $(\dim\Omega_i=\dim W^s (x))$\footnote{Here and throughout the article, $\dim$ denotes the topological dimension.}, where $x\in\Omega_i$, then the attractor (repeller) $\Omega_i$ is called expanding (contracting). It follows from \cite{2} (Theorem 2) that every expanding attractor (contracting repeller) $\Omega_i$ has the local structure of the direct product of the $r$-dimensional Euclidean space and the Cantor set, where $r$ is the topological dimension of $\Omega_i$.

Let $Diff (M^n)$ be the space of $C^1$ diffeomorphisms on $M^n$ endowed with the uniform $C^1$ topology \cite{22}. A diffeomorphism $f: M^n\rightarrow M^n$ is said to be structurally stable if there is a neighborhood $\mathcal U$ of the diffeomorphism $f$ in $Diff(M^n)$ such that every diffeomorphism $g\in\mathcal U$ is topologically conjugate to $f$. A diffeomorphism $f: M^n\rightarrow M^n$ is said to be $\Omega$-stable if there is a neighborhood $\mathcal U$ of the diffeomorphism $f$ in $Diff(M^n)$ such that for any $g\in\mathcal U$ restrictions $g|_{NW(g)}$ and $f|_{NW (f)}$ are topologically conjugate.

Let us introduce the relation $\prec$ for basic sets as follows: $\Omega_i\prec\Omega_j \Leftrightarrow W^s_{\Omega_i}\cap W^u_{\Omega_j}\not=\varnothing$. A $k$-cycle $(k\ge 1)$ is a collection of mutually disjoint basic sets $\Omega_0,\Omega_1,...,\Omega_k$ that satisfy the condition $\Omega_0\prec\Omega_1\prec...\prec\Omega_k\prec\Omega_0$. It follows from  \cite{19}, \cite{18} that the diffeomorphism $f: M^n\rightarrow M^n$ is $ \Omega$-stable if and only if it satisfies the axiom $A$ and has no cycles (for the formulation, see \cite{5}, Theorem 1.9.).

Let $M^2$ be a closed smooth orientable surface, $f: M^2\rightarrow M^2$ be a diffeomorphism satisfying the axiom $A$.

Let $\Lambda$ be a basic set of the diffeomorphism $f$. For $\sigma\in\{s, u\}$, let us put $\overline\sigma=s$ if $\sigma=u$, and $\overline\sigma=u$ if $\sigma=s$. A periodic point $p$ belonging to the set $\Lambda$ is called a boundary periodic point of type $\sigma$ ($\sigma$-boundary periodic point) if one of the connected components of the set $W^{\sigma}(p)\backslash p$ does not intersect with $\Lambda$, and both connected components of the set $W^{\overline\sigma}(p)\backslash p$ intersect with $\Lambda$. A periodic point $p$ belonging to the set $\Lambda$ is called a boundary periodic point of type $(s, u)$ if one of the connected components of each of the sets $W^s(p)\backslash p$, $W^u (p)\backslash p$ does not intersect with $\Lambda$. If the set $\Lambda$ is zero-dimensional or one-dimensional, then the set of boundary periodic points is non-empty and is finite (\cite{8}, \cite{9}).

If a basic set $\Lambda$ of the diffeomorphism $f$ is one-dimensional, then by virtue of \cite{2} it is an attractor or a repeller. If a one-dimensional basic set is an attractor (repeller), then it contains only $s$-boundary ($u$-boundary) periodic points.

Let $x$ be an arbitrary point belonging to a periodic component $\Lambda_i$ of a one-dimensional attractor (repeller) $\Lambda$. Then the set $W^u_x$ ($W^s_x$) belongs to the set $\Lambda_i$ and is dense in this set. Due to this fact and the fact that the closure of a connected set is connected, it follows that the periodic component $\Lambda_i$ is connected. Thus, each connected component of the one-dimensional attractor (repeller) $\Lambda$ coincides with one of its periodic components.

It is known (\cite{8}, \cite{10}) that for a one-dimensional attractor (repeller) $\Lambda$ accessible from inside boundary\footnote{Let $A$ be a subset of a topological space $X$. A point $y\in \partial A$ is called accessible from a point $x\in int$ $A$ if there exists a path $c: [0;1]\rightarrow X$ such that $c(0)=x, c (1)=y$ and $c (t)\in int$ $A$ for every $t\in (0;1)$. The union of all points accessible from the points of the set $int$ $A$ is called the boundary accessible from inside of the set $A$.} of the set $M^2\backslash\Lambda$ decays uniquely into a finite number of bunches. A bunch $b$ of the one-dimensional attractor $\Lambda$ is the union of the maximum number $h_b$ of the unstable manifolds $W^u_{p_1},..., W^u_{p_{h_b}}$ of the $s$-boundary periodic points $p_1,..., p_{h_b}$ of the set $\Lambda$ accessible from some (the same for all) point $x\in (M^2\backslash\Lambda)$. The number $h_b$ is called the degree of the bunch. Similarly, the concept of a bunch can be defined for a one-dimensional repeller.

In \cite{11}, for $A$-diffeomorphisms of compact surfaces (orientable and nonorientable), estimates of the maximum number of their one-dimensional basic sets are given, and the estimates are precise. It follows from \cite{11} that the maximum number of one-dimensional basic sets of the diffeomorphisms under consideration depends on the topological properties of the ambient surface (the genus of the surface and the number of connected components of the boundary) and on the geometric properties of one-dimensional basic sets (the number of bunches of degree one in one-dimensional basic sets).

In \cite{4}, a class $\mathbb G (M^2)$ of $A$-diffeomorphisms of closed orientable surfaces such that their non-wandering sets consist of one-dimensional basic sets is introduced, and necessary and sufficient conditions for the existence of a homeomorphism of the entire surface conjugating the restrictions of these diffeomorphisms on their non-wandering sets are found. The results of \cite{9} (Theorem 2.2. and Corollary to Theorem 3.2.) imply that a two-dimensional sphere and a two-dimensional torus do not admit diffeomorphisms from the class $\mathbb G (M^2)$. The first example of a diffeomorphism from the class $\mathbb G (M^2)$ was constructed in the work \cite{17}. Specifically, based on the $DA$-diffeomorphism of a two-dimensional torus (see \cite{20}, \cite{21}, \cite{1}, \cite{3}) and the diffeomorphism inverse to it, a diffeomorphism of a closed orientable surface of genus 2 (pretzel) such that its non-wandering set consists of a one-dimensional attractor and a one-dimensional repeller was constructed. It is proved in \cite{17} that the constructed diffeomorphism is not finitely $C^2$-stable (see the definition in \cite{17}), but is $\Omega$-stable. The primary mission of the given article is to research the interrelation between the dynamical properties of diffeomorphisms from the class $ \mathbb G(M^2)$ and the topology of the ambient surface $M^2$, as well as to research the stability of diffemorphisms from the class $\mathbb G (M^2)$.

Let $f:M^2\rightarrow M^2$ be a diffeomorphism from the class $\mathbb G (M^2)$ such that  its non-wandering set consists of $k_f$ periodic components $\Lambda_1,...,\Lambda_{k_f}$. It follows from the Lemma \ref{l123} that the non-wandering set of the diffeomorphism $f$ contains at least one attractor and at least one repeller, that is, $k_f\ge 2$. Let us denote by $m_{\Lambda_i}$ the number of bunches belonging to $\Lambda_i$, by $h_{\Lambda_i}$ the sum of the degrees of these bunches. Let us denote by $m_f$ the number of all bunches belonging to periodic components of the diffeomorphism $f$, by $h_f$ the sum of the degrees of these bunches. For the number $g\ge 0$, let us denote by $M^2_g$ a closed orientable surface of genus $g$.

\begin{theorem}\label{theorem1}
Let $f\in\mathbb G (M^2)$. Then the surface $M^2$ is homeomorphic to the connected sum:

$$M^2_{g_1} \#...\#M^2_{g_{k_f}}\#\underbrace{\mathbb T^2\#...\#\mathbb T^2}_{\text{$l_f$}},$$

where $g_i = 1+\frac{h_{\Lambda_i}}{4}-\frac{m_{\Lambda_i}}{2}$ $(i\in\{1,...,k_f\})$, $l_f = \frac{m_f}{2}-k_f+1$.
\end{theorem}

\begin{remark}
It follows from the Lemma \ref{l200} and the proof of the Theorem \ref{theorem1} that for every surface $M^2_{g_i}$, $i\in\{1,..., k_f\}$, there exists a compact submanifold $N_{\Lambda_i}\subset M^2$ of the genus $g_{\Lambda_i} = 1+\frac{h_{\Lambda_i}}{4}-\frac{m_{\Lambda_i}}{2}$ with $m_{\Lambda_i}$ boundary components which contains the periodic component $\Lambda_i$. Herewith, $N_{\Lambda_i}\cap N_{\Lambda_j}=\varnothing$ when $i\not= j$.
\end{remark}

\begin{corollary}
Let $f\in\mathbb G (M^2)$. Then the surface $M^2$ has the genus $g=1+\frac{h_f}{4}$.
\end{corollary}

\begin{figure}[h]
\center{\includegraphics[width=0.95\linewidth]{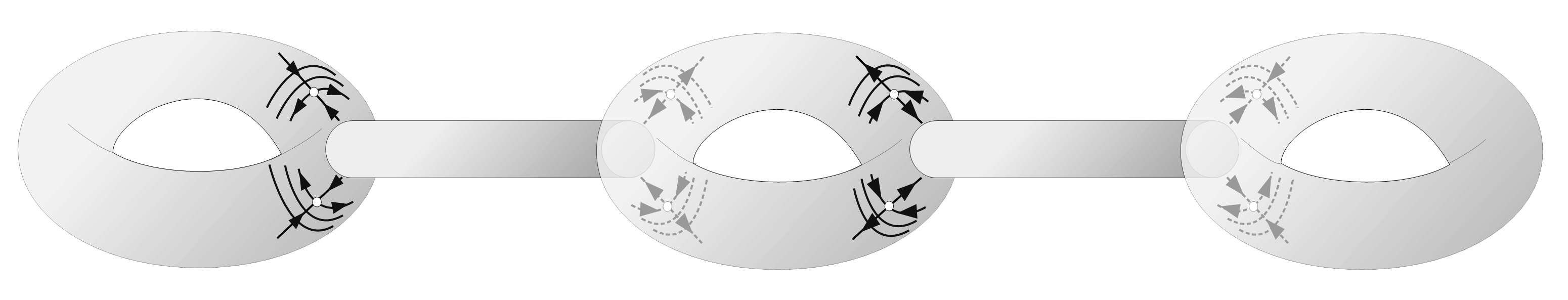} \\ a)}
\vfill
\center{\includegraphics[width=0.95\linewidth]{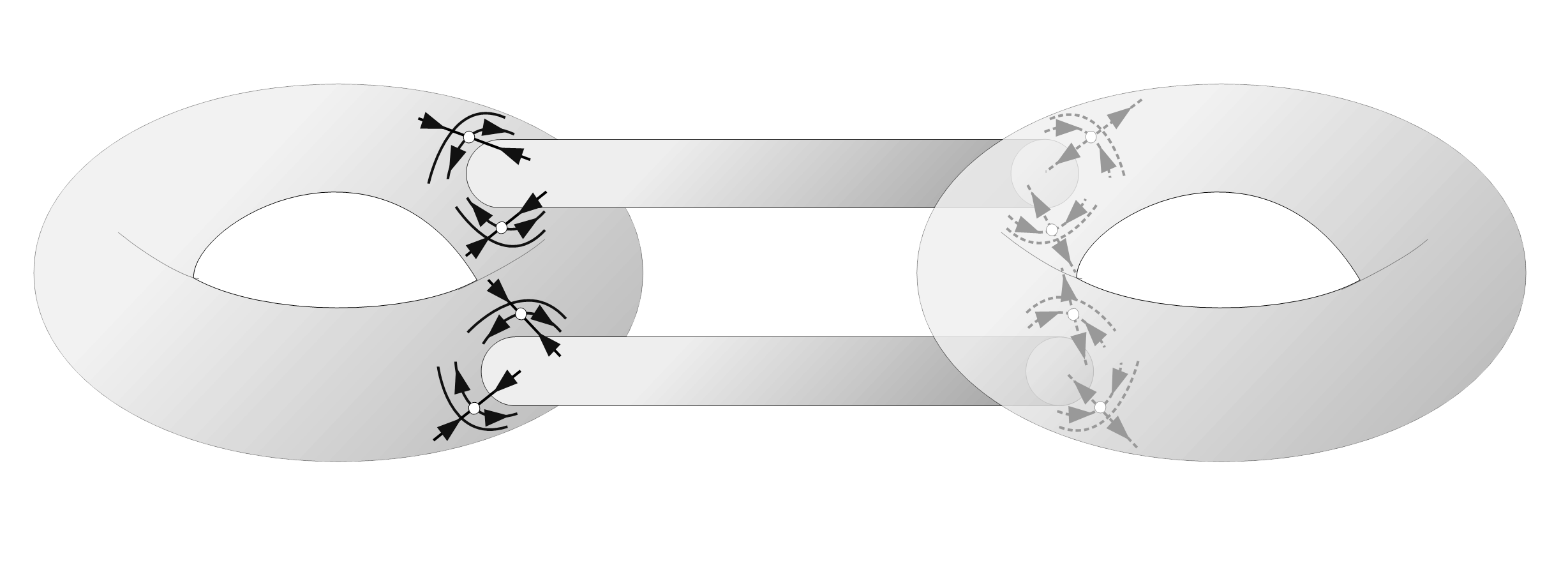} \\ b)}
\caption{Phase portrait of the diffeomorphism a) $f_1$; b) $f_2$.}
\label{ris20}
\end{figure}

Making use of the idea of constructing of the example from \cite{17}, one can construct examples of diffeomorphisms from the class $\mathbb G (M^2)$ on any closed orientable surface of genus $g\ge 2$. In the Figure \ref{ris20} a) it is shown a phase portrait of $A$-diffeomorphism $f_1$ of a closed orientable surface such that its non-wandering set consists of two one-dimensional attractors (each attractor has one bunch of degree two) and a one-dimensional repeller (which has two bunches of degree two). It follows from the Theorem \ref{theorem1} that $l_{f_1}=0$ and the ambient surface $M^2$ of the diffeomorphism $f_1$ is homeomorphic to the connected sum $M^2_{g_1} \#M^2_{g_2}\#M^2_{g_3}$, where $g_1=g_2=g_3=1$. In the Figure \ref{ris20} b) it is shown a phase portrait of $A$-diffeomorphism $f_2$ of a closed orientable surface such that its non-wandering set consists of a one-dimensional attractor (which has two bunches of degree two) and a one-dimensional repeller (which also has two bunches of degree two). It follows from the Theorem \ref{theorem1} that $l_{f_2}=1$ and the ambient surface $M^2$ of the diffeomorphism $f_2$ is homeomorphic to the connected sum $M^2_{g_1} \#M^2_{g_2}\#\mathbb T^2$, where $g_1=g_2=1$.

\begin{theorem}\label{theorem2}
Let $f\in\mathbb G(M^2)$. Then $f$ is $\Omega$-stable, but is not structurally stable.
\end{theorem}

\begin{ac}
This work was financially supported by the Russian Science Foundation (Project No. 21-11-00010), except for the proofs of Lemma \ref{l100} and Theorem \ref{theorem2}. The proof of Lemma \ref{l100} was obtained with the financial support from the Academic Fund Program at the HSE University in 2021– 2022 (grant № 21-04-004). The proof of Theorem \ref{theorem2} was obtained with the financial support from the Laboratory of Dynamical Systems and Applications NRU HSE, of the Ministry of science and higher education of the RF grant ag. No 075-15-2019-1931.
\end{ac}

\section{Auxiliary information and results}

\subsection{Surfaces with boundary}\label{2.1}

Let us recall that a compact two-dimensional manifold with a non-empty boundary is called a surface with boundary. The boundary of a compact surface is the union of finitely many mutually disjoint simple closed curves. Let us denote by $M^2 (Q)$ a compact surface with boundary, where $Q$ is the union of all simple closed curves belonging to the boundary of this surface. By gluing a finite number of closed two-dimensional disks (along their boundaries) to all the components of the boundary of the surface $M^2 (Q)$, one obtains a compact surface without boundary (a closed surface), which we denote by $M^2$. A surface $M^2 (Q$) is said to be orientable if the corresponding surface $ M^2$ without boundary is orientable. By definition, the genus of the surface $M^2 (Q)$ is equal to the genus of the surface $M^2$, and the Euler characteristic of $M^2 (Q)$ is equal to the difference between the Euler characteristic of $M^2$ and the number of curves in the set $Q$.

\subsection{Saddle singularities}\label{2.3}

Let $k\in\mathbb N$ and $k\not=2$. The foliation $W_k$ on $\mathbb R^2$ with the standard saddle singularity at the point $O$ (coordinate origin) and $k$ separatrices is the image of the horizontal lines $\{Im\, w=c, c\in \mathbb R\}$ under the map $w=z^{\frac{k}{2}}$ in the case of odd $k$ and under the map $w^2=z^k$ in the case of even $k$. For $k=2$ all the leaves of the foliation $W_2 $ are straight lines $y=c$, but the axis $Ox$ is artificially split into three parts: the origin and two the half-axes, the latter called the separatrices.

Let $M^2$ be a closed surface, $\mathcal F$ be a foliation on the surface $M^2$. The foliation $\mathcal F$ is said to be a foliation with saddle singularities if the set $\mathcal S$ of the singularities of the foliation $\mathcal F$ consists of a finite number of points and for any point $s\in\mathcal S$ there is a neighborhood $U_s\subset M^2$, the homeomorphism $\psi_s:U_s\rightarrow \mathbb R^2$ and the number $k_s\in\mathbb N$ such that $\psi_s(s)=O$ and $\psi_s(\mathcal F\cap U_s)=W_{k_s}\backslash O$. The point $s$ is called the saddle singularity with $k_s$ separatrices. Index $I (s)$ of each saddle singularity $s\in\mathcal S$ can be calculated via the number of separatrices $k_s$ by the formula:

\begin{equation}\label{form1}
I(s)=1-\frac{k_s}{2}. 
\end{equation}

Let $\chi(M^2)$ be the Euler characteristic of the surface $M^2$. The next formula follows from the Poincar\'e-Hopf theorem:

\begin{equation}\label{form2}
\chi(M^2)=\sum\limits_{s\in\mathcal S} I(s).
\end{equation}

\subsection{Auxiliary lemmas}\label{2.2}

Let $M^2$ be a closed smooth orientable surface, $f: M^2\rightarrow M^2$ be an $A$-diffeomorphism such that its non-wandering set contains a one-dimensional attractor (repeller). Let $\Lambda$ be a periodic component of this attractor (repeller), $b_1,..., b_{m_{\Lambda}}$ be the bunches belonging to $\Lambda$ ($m_{\Lambda}$ bunches in total), $h_{\Lambda}$ be the sum of the degrees of these bunches.

The proof of the following lemma uses the ideas from \cite{15}, \cite{6}, \cite{7}, as well as the proof scheme from \cite{5} (Theorem 9.6.).

\begin{lemma}\label{l200}
For the periodic component $\Lambda$ of a one-dimensional attractor (repeller) of the diffeomorphism $f: M^2\rightarrow M^2$, there are a submanifold $N_{\Lambda}$ and a natural number $n$ with the following properties:

\begin{enumerate}
\item $N_{\Lambda}$ is a trapping neighborhood of the set $\Lambda$ with respect to the diffeomorphism $f^n$;

\item $N_{\Lambda}$ is a compact orientable surface of the genus $g_{\Lambda} = 1+\frac{h_{\Lambda}}{4}-\frac{m_{\Lambda}}{2}$ with $m_{\Lambda}$ boundary components.
\end{enumerate}
\end{lemma}

\begin{figure}[h!]
\center{\includegraphics[scale=0.4]{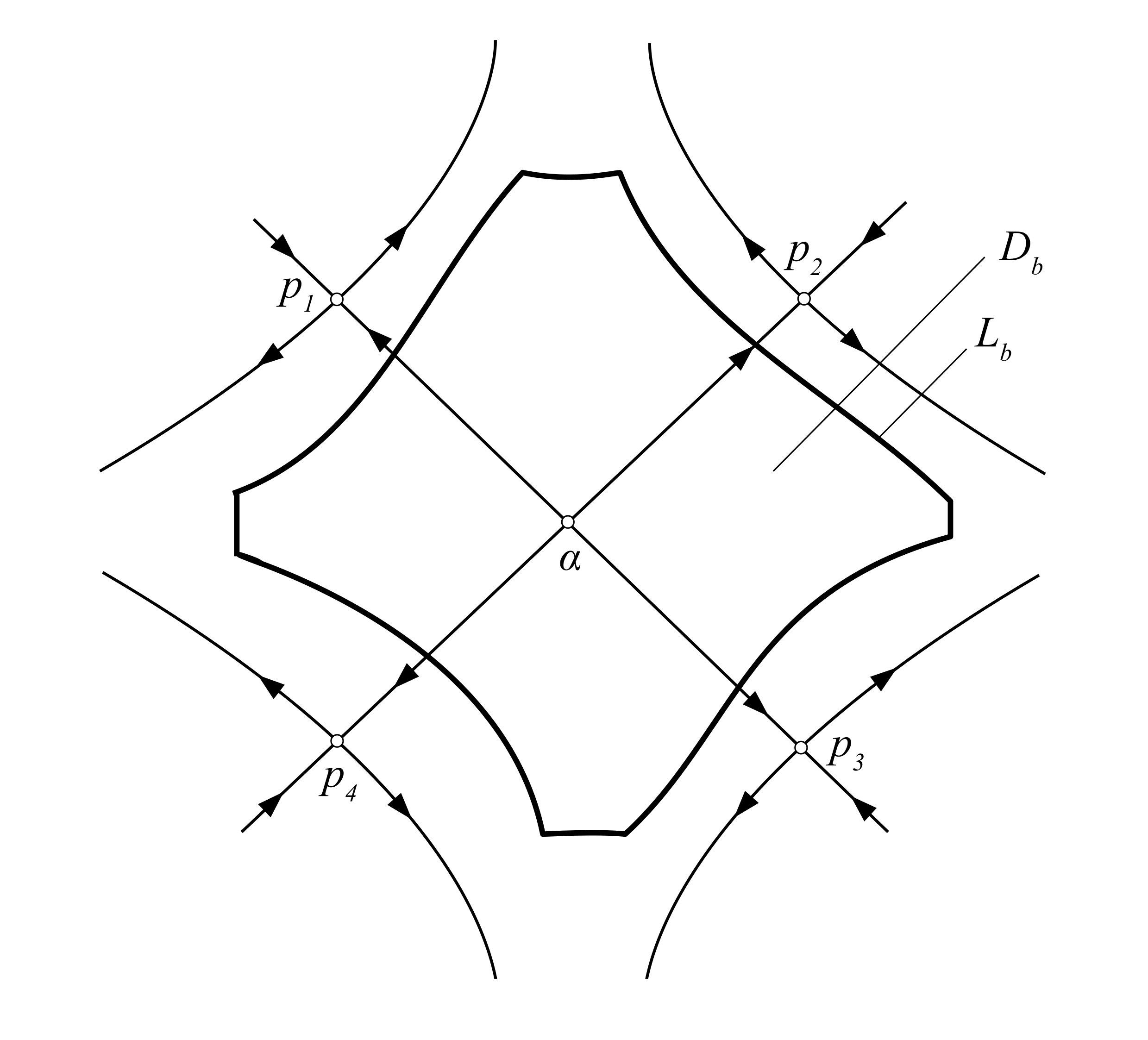}}
\caption{Construction of the surface $M_{\Lambda}$.}
\label{r101}
\end{figure}

\begin{demo}

For definiteness, we will assume that $\Lambda$ is a periodic component of the attractor of the diffeomorphism $f$ (if $\Lambda$ is a periodic component of the repeller, it is sufficient to consider the diffeomorphism $f^{-1}$).

The finiteness of the number of periodic components of a basic set and the finiteness of the set of boundary periodic points of a one-dimensional attractor imply that there exists a number $n\in\mathbb N$ such that $f^n (\Lambda)=\Lambda$ and all boundary periodic points of the set $\Lambda$ are fixed with respect to the diffeomorphism $f^n$.

Let us denote by $b$ an arbitrary bunch belonging to the set $\Lambda$, by $h_b$ the degree of this bunch. It follows from the definition of the bunch that $b=W^u_{p_1}\cup...\cup W^u_{p_{h_b}}$, where $p_j, j\in\{1,..., h_b\}$, is $s$-boundary periodic point of the set $\Lambda$. By virtue of \cite{8} (Lemma 3.3), there exists a sequence of points $x_1,..., x_{2h_{b}}$ such that:

\begin{enumerate}
\item $x_{2j-1}, x_{2j}$ belong to different connected components of the set $W^u_{p_j}\backslash p_j$;

\item $x_{2j+1}\in W^s_{x_{2j}}$ (we assume $x_{2h_b+1}=x_1$);

\item $(x_{2j},x_{2j+1})^s\cap\Lambda=\varnothing, j=1,...,h_b$.
\end{enumerate}

For each $j\in\{1,..., h_b\}$, let us choose a pair of points $\tilde x_{2j-1}, \tilde x_{2j}$, and a simple curve $l_j$ with the boundary points $\tilde x_{2j-1}, \tilde x_{2j}$ such that:

\begin{enumerate}
\item $(\tilde x_{2j},\tilde x_{2j+1})^s \subset(x_{2j},x_{2j+1})^s$ $(x_{2h_b+1}=x_1)$;

\item the curve $l_j$ transversally intersects with the stable manifold of any point belonging to the arc $(x_{2j-1}, x_{2j})^u$ at exactly one point;

\item $L_b=\bigcup\limits_{j=1}^{h_b}[l_j\cup(\tilde x_{2j},\tilde x_{2j+1})^s]$ is a simple closed piecewise smooth curve and the set $L_{\Lambda} = \bigcup\limits_{t\in\{1,..., m_{\Lambda}\}} L_{b_t}$ has the properties:

\begin{enumerate}
\item $f^n(L_{\Lambda})\cap L_{\Lambda}=\varnothing$;

\item for every curve $L_{b_t}, t\in\{1,..., m_{\Lambda}\}$, there exists a curve from the set $f^n (L_{\Lambda})$ such that these two curves are the boundary of the two-dimensional annulus $K_{b_t}$;

\item the annuli $\{K_{b_t},t\in\{1,...,m_{\Lambda}\}\}$ are pairwise disjoint (see Figure \ref{r100}).
\end{enumerate}
\end{enumerate}

For an arbitrary bunch $b$, we will call the curve $L_b$ the characteristic curve of the bunch $b$.

\begin{figure}[h!]
\center{\includegraphics[scale=0.4]{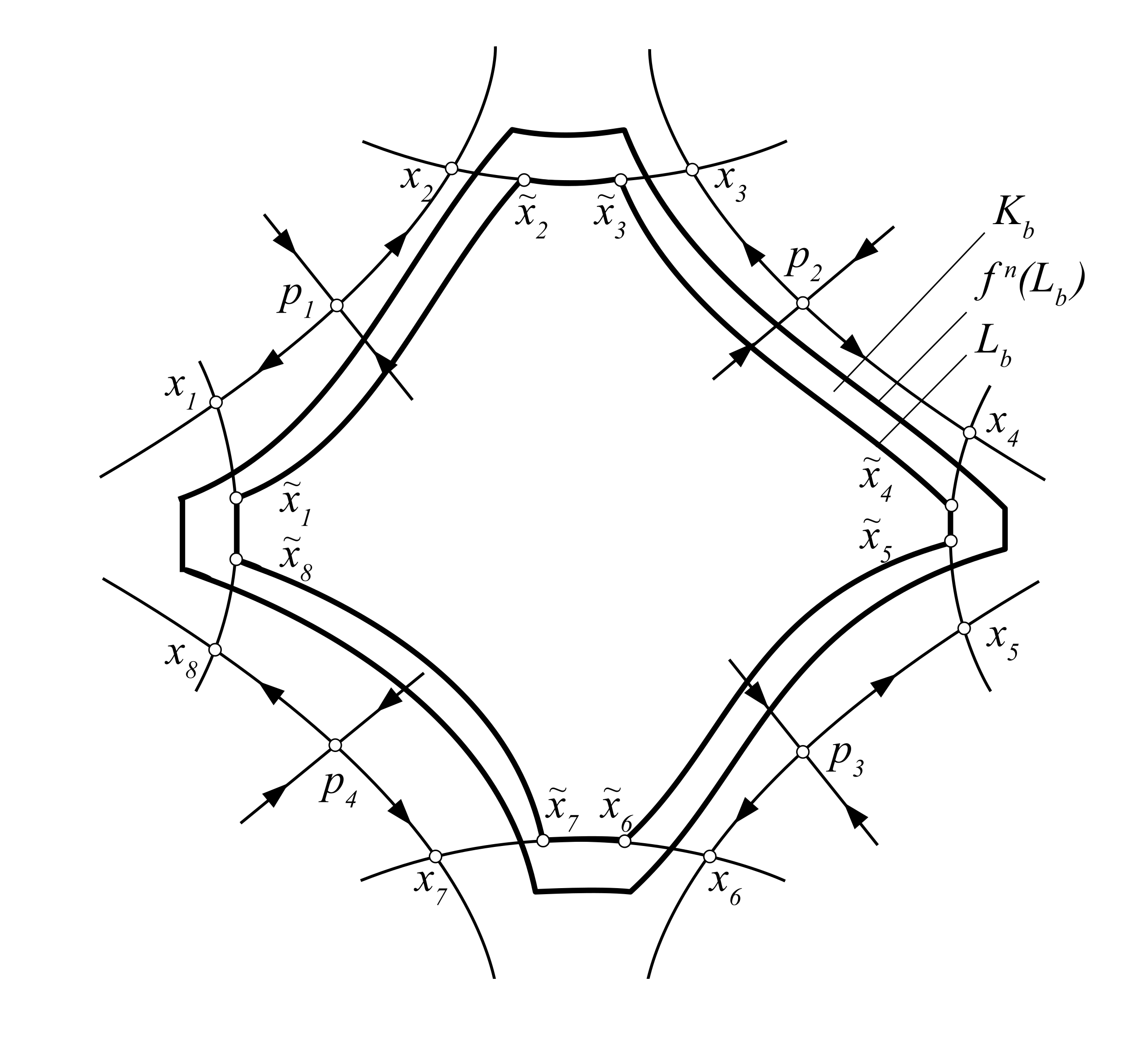}}
\caption{Construction of the characteristic curve for the bunch of degree 4.}
\label{r100}
\end{figure}

Let us put $N_{\Lambda} = \Lambda\cup\bigcup\limits_{k\ge 1} f^{kn}(\bigcup\limits_{t\in\{1,..., m_{\Lambda}\}}K_{b_t})$. By construction, the annuli $\{K_{b_t}, t\in\{1,..., m_{\Lambda}\}\}$ consist of wandering points of the diffeomorphism $f^n$ and $N_{\Lambda}$ is a compact orientable surface with non-empty boundary (consisting of $m_{\Lambda}$ components) such that $f^n(N_{\Lambda})\subset int\, N_{\Lambda}$ and $\Lambda=\bigcap\limits_{k\ge 0} f^{kn}(N_{\Lambda})$. Thus, $N_{\Lambda}$ is the trapping neighborhood of the set $\Lambda$ with respect to the diffeomorphism $f^n$. Since the set $\Lambda$ is connected and $N_{\Lambda}$ is its trapping neighborhood, then $N_{\Lambda}$ is connected.

We will prove that the genus $g_{\Lambda}$ of the surface $N_{\Lambda}$ is equal to $1+\frac{h_{\Lambda}}{4}-\frac{m_{\Lambda}}{2}$. Remove the set $int(\bigcup\limits_{t\in \{1,..., m_{\Lambda}\}}K_{b_t})$ from the surface $M^2$. As a result, the surface $M^2$ decomposes into a finite number of connected components, one of which is the set $N_{\Lambda}$. In this case, the set $\bigcup\limits_{t\in\{1,..., m_{\Lambda}\}} f^n(L_{b_t})$ is the boundary of the set $N_{\Lambda}$. To each curve $f^n (L_{b_t})$ $(t\in\{1,..., m_{\Lambda}\})$ let us glue a closed two-dimensional disk $D_{b_t}$ (along its boundary) and denote the obtained manifold by $M_{\Lambda}$. Let us construct a homeomorphism $F: M_{\Lambda}\rightarrow M_{\Lambda}$ such that $F|_{N_{\Lambda}}=f|_{N_{\Lambda}}$ and the non-wandering set of $F|_{D_{b_t}}$ (for all $t\in\{1,..., m_{\Lambda}\}$) consists of exactly one hyperbolic periodic source point $\alpha_{b_t}$. By construction, $\alpha_{b_t}$ belongs to the closure $W^s_{p_j}$ for each $j\in\{1,..., h_{b_t}\}$ (see Figure \ref{r101}).

Let us put $S_{\Lambda}=\bigcup\limits_{t\in \{1,..., m_{\Lambda}\}}\alpha_{b_t}$. The surface $M_{\Lambda}$ admits a foliation

$$\mathcal F_{M_{\Lambda}}=\{W^s_x,x\in(\Lambda\cup S_{\Lambda})\},$$

which has $m_{\Lambda}$ singularities (points $\alpha_{b_t}, t\in\{1,...,m_{\Lambda}\}$), and all these singularities are saddle. The formula \ref{form1} (see section \ref{2.3}) implies that the index $I (\alpha_{b_t})$ of each saddle singularity $\alpha_{b_t}$ is equal to $(1-\frac{h_{b_t}}{2})$. From here and from the formula \ref{form2} (see section \ref{2.3}) one gets:

\begin{equation}\label{form3}
\chi(M_{\Lambda})=\sum\limits_{t\in \{1,...,m_{\Lambda}\}}I(\alpha_{b_t})=m_{\Lambda}-\frac{h_{\Lambda}}{2},  
\end{equation}

\text{where $\chi(M_{\Lambda})$ is the Euler characteristic of the surface $M_{\Lambda}$}.

Since $M_{\Lambda}$ is closed orientable surface, its genus $g_{\Lambda}$ is related to the Euler characteristic $\chi (M_{\Lambda})$ by the following formula: $\chi(M_{\Lambda})=2-2g_{\Lambda}$. This fact and the formula (\ref{form3}) imply that the genus of the surface $M_{\Lambda}$ is calculated by the formula $g_{\Lambda}=1+\frac{h_{\Lambda}}{4}-\frac{m_{\Lambda}}{2}$.

It follows from the construction of the surface $M_{\Lambda}$ that $N_{\Lambda}=M_{\Lambda}\backslash(\bigcup\limits_{t\in \{1,..., m_{\Lambda}\}} int \, D_{b_t})$. Hence, the surface $N_{\Lambda}$ has the same genus as $M_{\Lambda}$.

\end{demo}

Let $M^2$ be a closed orientable surface and $\varphi: S\times[0;1]\rightarrow M^2$ be an embedding\footnote{The map $\varphi: X \rightarrow Y$, where $X, Y$ are topological spaces, is said to be an embedding if $\varphi: X\rightarrow \varphi (X)\subset Y$ is a homeomorphism, where $\varphi (X)$ carries the subspace topology inherited from $Y$.}, where $S=\{(x_1, x_2)\in\mathbb R^2| \; x_1^2+x_2^2=1\}$ is a unit circle. Let us denote by $U$ the image of the set $S\times (0;1)$ with respect to the map $\varphi$.

The proof of the following lemma uses the ideas of the proofs of Lemmas 6 and 7 \cite{12}.

\begin{lemma}\label{l1}
If the manifold $M^2\backslash U$ is connected, then there exists a closed orientable surface $M^2_1$ such that $M^2$ is homeomorphic to the connected sum $M^2_1\#\mathbb T^2$.
\end{lemma}

\begin{demo}
Let us denote by $S_1$ and $S_2$ the components of the boundary $\partial U$ of the set $U$. Note that $S_1$ and $S_2$ are homeomorphic to a circle. Let us denote by $\tilde M$ the manifold obtained as follows: $\tilde M = (M^2\backslash U)\cup(B_1\cup B_2)$, where $B_1$ and $B_2$ are disjoint closed two-dimensional disks such that $\partial B_1=S_1$, $\partial B_2=S_2$. Since the manifold $M^2\backslash U$ is connected, then the manifold $\tilde M$ is connected. It follows from the connection of $\tilde M$ that there is a closed subset $\tilde D\subset\tilde M$ such that it contains $B_1\cup B_2$ and is homeomorphic to a two-dimensional disk. Let us put $Q = \tilde D\backslash(int$ $B_1\cup$ $ int$ $B_2)$. The boundary of the set $Q$ consists of three connected components: $S_0,S_1, S_2$, each of which is homeomorphic to a circle. By gluing a closed two-dimensional disk (along its boundary) to the set $Q$ along the component $S_0$, one obtains a set homeomorphic to $S^2\backslash(int$ $B_1\cup$ $ int$ $B_2)$, where $S^2$ is a two-dimensional sphere.

To get  from $M^2\backslash U$ the original surface $M^2$, it is necessary to glue $S_1$ and $S_2$ by means of an orientation-reversing homeomorphism. At the same time, the manifold obtained from $S^2\backslash(int$ $B_1\cup$ $int$ $B_2)$ after this gluing is homeomorphic to a two-dimensional torus $\mathbb T^2$. Thus, the original surface $M^2$ is homeomorphic to the connected sum $M^2_1\# \mathbb T^2$, where $M^2_1$ is a closed orientable surface.
\end{demo}

\begin{lemma}\label{l123}
Let $f\in\mathbb G(M^2)$. Then its non-wandering set contains at least one attractor and at least one repeller.
\end{lemma}

\begin{demo}
Assume the opposite. Let the non-wandering set of the diffeomorphism $f$ consist of one-dimensional attractors $\Lambda_1,...,\Lambda_{k_f}$ (if it consists of one-dimensional repellers, then it is sufficient to consider the diffeomorphism $f^{-1}$). According to \cite{2} (Theorem 2), the set $\Lambda_i$ $(i=\overline{1,...,k_f})$ has the local structure of the direct product of the interval and the Cantor set. Thus, every set $\Lambda_i$ is nowhere dense. The properties of basic set that is an attractor imply that $\Lambda_i=W^u(\Lambda_i), i=\overline{1,...,k_f}$. It follows from \cite{1} (Corollary 6.3) and aforesaid that $M^2=\bigcup\limits_{i\in\{1,...,k_f\}} W^u(\Lambda_i) = \bigcup\limits_{i\in\{1,...,k_f\}}\Lambda_i$. That contradicts Baire category theorem which states that a non-empty complete metric space cannot be represented as a countable union of nowhere dense subsets.
\end{demo}

\section{Proof of the main results}

Throughout this section, $f:M^2\rightarrow M^2$ is a diffeomorphism from the class $\mathbb G (M^2)$. Let us denote the periodic components of this diffeomorphism by $\Lambda_i$ $(i\in\{1,..., k_f\})$, the set $M^2\backslash \bigcup\limits_{i=1}^{k_f} \Lambda_i$ by $V$. One can choose the number $m\in\mathbb N$ such that $f^m(\Lambda_i)=\Lambda_i$ for all $i\in\{1,..., k_f\}$ and all boundary periodic points of all one-dimensional basic sets are fixed with respect to the diffeomorphism $f^m$. Therefore, without loss of generality, throughout this section we will assume that every basic set has a unique periodic component and all boundary periodic points of all one-dimensional basic sets are fixed with respect to the diffeomorphism $f$. We will call each set $\Lambda_i$ a basiс set.

\begin{lemma}\label{l100}
The set $V$ consists of a finite number of mutually disjoint open connected sets $V_1,..., V_l$ such that the boundary accessible from inside of each such set consists of two bunches, one of which belongs to some attractor, and the other belongs to some repeller of the diffeomorphism $f$.
\end{lemma}

\begin{demo}
It follows from \cite{16} (see proposition 2 and corollary to it) that if $A\subset M$ is an open subset of a closed manifold $M$ (it is possible that $A$ coincides with $M$), $B\subset A$ is a closed subset of $A$, then each connected component of the complement $A\backslash B$ is an open set.

Let us consider the set $M^2\backslash\Lambda_i$, where $\Lambda_i$ is one of the basic sets of the diffeomorphism $f$. The aforesaid and the closedness of the set $\Lambda_i$ imply that each connected component of the complement $M^2\backslash \Lambda_i$ is an open set. From this fact and from the fact that the boundary accessible from inside of the set $M^2\backslash\Lambda_i$ decays uniquely into a finite number of bunches (herewith, each bunch is the boundary accessible from inside of exactly one connected set) (see \cite{8}, \cite{10}), it follows that the set $M^2\backslash\Lambda_i$ consists of a finite number of mutually disjoint open connected sets. Let $\Lambda_j$ be the basic set of the diffeomorphism $f$ such that $\Lambda_j\not =\Lambda_i$. Applying the same reasoning and considering that $\Lambda_i\cap\Lambda_j=\varnothing$, one gets that $M^2\backslash(\Lambda_i\cup\Lambda_j)$ consists of a finite number of mutually disjoint open connected sets. Consistently applying these arguments $k_f$ times (for all basic sets of the diffeomorphism $f$), one obtains that the set $V=M^2\backslash \bigcup\limits_{i=1}^{k_f} \Lambda_i$ consists of a finite number of mutually disjoint open connected sets $V_1,...,V_l$.

Let $x$ be an arbitrary point belonging to the set $V_i$ $(i\in\{1,...,l\})$. According to \cite{1} (Corollary 6.3), the point $x$ lies on the stable manifold of some basic set $\Lambda^a_i$ and on the unstable manifold of some basic set $\Lambda^r_i$. Since $x\in V$, the set $NW(f)$ consists of one-dimensional attractors and one-dimensional repellers and stable (unstable) manifold of the repeller (attractor) belongs to it, then the set $\Lambda^a_i$ is an attractor, and the set $\Lambda^r_i$ is a repeller of the diffeomorphism $f$. The connection of the set $V_i$ implies that all points belonging to this set lie on the stable manifold of the attractor $\Lambda^a_i$ and on the unstable manifold of the repeller $\Lambda^r_i$. From this fact, from the connection of the set $V_i$ and from the fact that the boundary accessible from inside of the sets $M^2\backslash\Lambda^a_i$ and $M^2\backslash\Lambda^r_i$ decays uniquely into a finite number of bunches, it follows that the boundary accessible from inside of the set $V_i$ consists of two bunches, one of which belongs to the attractor $\Lambda^a_i$, and the other belongs to the repeller $\Lambda^r_i$.
\end{demo}

The proof of the Lemma \ref{l100} implies the following Corollary.

\begin{corollary}\label{cor1}
The number of bunches of all attractors of the diffeomorphism $f$ is equal to the number of bunches of all its repellers.
\end{corollary}

\begin{demo1}

Let $\Lambda^a$ be a one-dimensional attractor of the diffeomorphism $f$, $b^a$ be one of its bunches, $L_{b^a}$ be the characteristic curve of the bunch $b^a$ (see the proof of the Lemma \ref{l200}). Let $S=\{(x_1, x_2)\in\mathbb R^2| \; x_1^2+x_2^2=1\}$ be a unit circle, $\varphi: S\times[0;1]\rightarrow M^2$ be an embedding such that $\varphi (S\times\{0\})=L_{b^a}, \varphi (S\times\{1\})=f(L_{b^a})$. Let us denote by $U$ the image of the set $S\times (0;1)$ with respect to the map $\varphi$. It follows from the proof of the Lemma \ref{l100} that there exists $i\in\{1,..., l\}$ such that the curve $L_{b^a}$ belongs to the set $V_i$. Hence, $U\subset V_i$. Remove the set $U$ from the surface $M^2$. The boundary of $M^2 \backslash U$ consists of two connected components, each of which is homeomorphic to a circle. Let us glue closed two-dimensional disks $B_1$ and $B_2$ (along their boundaries) to these components and denote the resulting surface by $\tilde M^2$. There are two possible cases:

\begin{enumerate}
\item the surface $\tilde M^2$ is connected. Then, according to Lemma \ref{l1}, there exists a closed orientable surface $P$ such that the surface $M^2$ is homeomorphic to the connected sum $P\# \mathbb T^2$;

\item the surface $\tilde M^2$ is disconnected. Then the surface $M^2$ is homeomorphic to the connected sum $P_1\# P_2$, where $P_1$ and $P_2$ are closed orientable surfaces.
\end{enumerate}

It follows from the Lemma \ref{l100} that the boundary accessible from inside of the set $V_i$ consists of the bunch $b^a$ of the attractor $\Lambda^a$ and a bunch $b^r$ of some repeller $\Lambda^r$ of the diffeomorphism $f$. We denote by $p_1,...,p_{h_{b^a}}$ the boundary periodic points belonging to the bunch $b^a$, and by $q_1,..., q_{h_{b^r}}$ the boundary periodic points belonging to the bunch $b^r$.

Let us define a homeomorphism $F: \tilde M^2\rightarrow \tilde M^2$ such that:

\begin{enumerate}
\item $F|_{\tilde M^2\backslash (int \, B_1\cup \, int \, B_2)}$ = $f|_{\tilde M^2\backslash (int \, B_1\cup \, int \, B_2)}$;

\item the non-wandering set of $F|_{B_1}$ consists of exactly one hyperbolic fixed source point $\alpha$ (by construction, this point belongs to the closure $W^s_{p_j}$ for each $j\in\{1,..., h_{b^a}\}$);

\item the non-wandering set of $F|_{B_2}$ consists of exactly one hyperbolic fixed sink point $\omega$ (by construction, this point belongs to the closure $W^u_{q_j}$ for each $j\in\{1,..., h_{b^r}\}$).
\end{enumerate}

Let us consistently perform the procedure described above for all bunches belonging to attractors of the diffeomorphism $f$. As a result, one gets a disconnected manifold, which is the union of $k_f$ closed orientable surfaces. It follows from the proof of Lemma \ref{l200} that every such surface has genus $g_i = 1+\frac{h_{\Lambda_i}}{4}-\frac{m_{\Lambda_i}}{2}$ $(i\in\{1,..., k_f\})$ (see the notation in the condition of the Theorem \ref{theorem1}). Since the number of all bunches belonging to attractors of the diffeomorphism $f$ is equal to $ \frac{m_f}{2}$ (see Corollary \ref{cor1}), then the procedure described above is performed $\frac{m_f}{2}$ times. Among them, there are $k_f-1$ steps, as a result of which the manifold splits into two disconnected manifolds, and $l_f = \frac{m_f}{2}-k_f+1$ steps, as a result of which the manifold remains connected. Thus, one obtains that the original surface $M^2$ is homeomorphic to the connected sum:

$$M^2_{g_1} \#...\#M^2_{g_{k_f}}\#\underbrace{\mathbb T^2\#...\#\mathbb T^2}_{\text{$l_f$}},$$

where $g_i = 1+\frac{h_{\Lambda_i}}{4}-\frac{m_{\Lambda_i}}{2}$ $(i\in\{1,...,k_f\})$, $l_f = \frac{m_f}{2}-k_f+1$.
\end{demo1}

\begin{demo3}
Since the surface $M^2$ is homeomorphic to the connected sum of $k_f$ closed orientable surfaces of the genus $g_i$ $(i\in\{1,..., k_f\})$ and $l_f$ two-dimensional tori, then the genus $g$ of the surface $M^2$ is calculated by the following formula:

$$g = \sum\limits_{i=1}^{k_f}(1+\frac{h_{\Lambda_i}}{4}-\frac{m_{\Lambda_i}}{2})+\frac{m_f}{2}-k_f+1=k_f+\frac{h_f}{4}-\frac{m_f}{2}+\frac{m_f}{2}-k_f+1=1+\frac{h_f}{4}.$$

\end{demo3}

\begin{demo2}
In \cite{6} (Theorem 1), it is proved that if the non-wandering set of a structurally stable diffeomorphism of a closed smooth orientable surface contains a one-dimensional attractor (repeller), then it contains a source (sink) periodic point. This fact and the fact that the non-wandering set of the diffeomorphism $f$ consists of one-dimensional attractors and one-dimensional repellers entail that the diffeomorphism $f$ is not structurally stable.

We will prove that the diffeomorphism $f$ is $\Omega$-stable. Let $\Lambda^a$ be an arbitrary one-dimensional attractor of the diffeomorphism $f$. Unstable manifold of this attractor coincides with it, and stable manifold of this attractor, by virtue of \cite{1} (Corollary 6.3), intersects with the unstable manifolds of a finite number of repellers $\Lambda^{r_1},...,\Lambda^{r_l}$. Herewith, the stable manifold of each of the repellers $\Lambda^{r_1},...,\Lambda^{r_l}$ coincides with it. Conducting similar reasoning for an arbitrary repeller of the diffeomorphism $f$, one obtains that the diffeomorphism $f$ has no cycles. Hence, according to \cite{5} (Theorem 1.9.), the diffeomorphism $f$ is $\Omega$-stable.
\end{demo2}

\renewcommand\refname{References}

\end{document}